\documentclass[10pt,twoside]{classecras-e}
\usepackage{amssymb,amsbsy,amsmath,amsfonts,amssymb,amscd}
\usepackage{latexsym,euscript}
\usepackage[english,francais]{babel}
\usepackage{times}
\input utile.def

\newtheorem{theoremempty}{Theorem}

\newtheorem{maintheorem}[theoremempty]{Main Theorem}
\newtheorem{maintheoremrestated}[theoremempty]{Main Theorem Restated}
\newtheorem{Theorem}{Theorem}
\newtheorem{Conj}{Conjecture}

 \def\RR{{\mathbb R}}

   \def\cN{{\cal N}}

\def\dim{\operatorname{dim}}

\ComParit{S0764-4442}
\PIT{FLA}
\PXMA{????}
\Add{?}
\Volume{332}
\Year{2001}
\FirstPage{1}
\LastPage{??}
\AuteurCourant{S. Crovisier and  D. Yang}
\TitreCourant{Singular hyperbolic three-dimensional vector fields} 
\Journal
\Rubrique{Syst\`emes dynamiques}{Dynamical systems}
\PresentePar{First name}{NAME}
\Recu{jour mois ann\'ee}{apr\`es r\'evision}{jour mois ann\'ee}
\begin{document}
\selectlanguage{english}
\title{On the density of singular hyperbolic three-dimensional vector fields: a conjecture of Palis
}
\author{%
Sylvain Crovisier~$^{\text{a}}$
\footnote{S.Crovisier was partially supported by \emph{Balzan Research Project} of J.Palis.},\ \
Dawei Yang~$^{\text{b}}\footnote{D. Yang was partially supported by NSFC 11271152
and ANR project \emph{DynNonHyp} BLAN08-2 313375}$
}
\address{%
\begin{itemize}\labelsep=2mm\leftskip=-5mm
\item[$^{\text{a}}$]
CNRS - Laboratoire de Math\'ematiques d'Orsay,
Universit\'e Paris-Sud 11, Orsay 91405, France\\
E-mail: Sylvain.Crovisier@math.u-psud.fr
\item[$^{\text{b}}$]
School of Mathematical Sciences, Soochow University, Suzhou, 215006, P.R. China\\
E-mail: yangdw1981@gmail.com, yangdw@suda.edu.cn
\end{itemize}
}
\maketitle
\thispagestyle{empty}
\begin{Abstract}{%
In this note we announce a result for vector fields on three-dimensional manifolds:
those who are singular hyperbolic or exhibit a homoclinic tangency
form a dense subset of the space of $C^1$-vector fields.
This answers a conjecture by Palis. The argument uses an extension
for local fibered flows of Ma\~n\'e and Pujals-Sambarino's theorems about the  uniform contraction
of one-dimensional dominated bundles.
}\end{Abstract}
\selectlanguage{french}
\begin{Ftitle}{%
Sur la densit\'e de l'hyperbolicit\'e singuli\`ere pour les champs de vecteurs en dimension trois~:
une conjecture de Palis
}\end{Ftitle}
\begin{Resume}{%
Dans cette note, nous annon\c{c}ons un r\'esultat portant sur les champs de vecteurs
des vari\'et\'es de dimension $3$~:
ceux qui v\'erifient l'hyperbolicit\'e singuli\`ere ou qui poss\`edent une tangence homocline
forment un sous-ensemble dense de l'espace des champs de vecteurs $C^1$.
Ceci r\'epond \`a une conjecture de Palis.
La d\'emonstration utilise une g\'en\'eralisation pour les flots fibr\'es locaux
des th\'eor\`emes de Ma\~n\'e et Pujals-Sambarino traitant de la contraction uniforme
de fibr\'es unidimensionnels domin\'es.
}\end{Resume}

\par\medskip\centerline{\rule{2cm}{0.2mm}}\medskip
\setcounter{section}{0}
\selectlanguage{english}

\section{Introduction}

We are concerned with dynamics which are typical in the space of dynamical systems. Hyperbolic systems are natural candidates since they form an open set. However there are obstructions for hyperbolicity, such as homoclinic tangencies which produces rich behavior as Newhouse phenomena \cite{New79}. Palis \cite{Pal91,Pal00} conjectured that on surfaces,
every diffeomorphism can be accumulated either by hyperbolic ones or by diffeomorphisms with  a homoclinic tangency. Pujals and Sambarino \cite{PuS00} managed to prove the conjecture in the $C^1$ topology. For higher dimensions, other homoclinic bifurcations have to be included in Palis conjecture, such as heterodimensional cycles.
For the case of vector fields, new kinds of bifurcations, associated to singularities,
and called singular cycles, have to be introduced. A version of Palis conjecture for vector fields can be formulated as:

\begin{Conj} \emph{(Palis \cite{Pal05}) -- }
Every vector field can be accumulated either by hyperbolic vector fields or by ones with a homoclinic bifurcation or with a singular cycle.
\end{Conj}
\medskip

Vector fields admit robustly non-hyperbolic transitive attractors with singularities \cite{ABS77,Guc76}. Morales, Pacifico and Pujals \cite{MPP04} defined the \emph{singular hyperbolicity} to characterize robust attractors with singularities in dimension $3$.
Palis gave a stronger version of his conjecture in dimension $3$
(see also \cite{ArR03,MoP03}):

\begin{Conj} \emph{(Palis \cite{Pal05}) -- }\label{conj2}
Every vector field on a three dimensional manifold can be accumulated either by singular hyperbolic vector fields or by ones with a homoclinic tangency (associated to a non-singular periodic orbit).
\end{Conj}
\medskip

Arroyo and Rodriguez-Hertz \cite{ArR03} proved the first conjecture for three-dimensional manifolds in the $C^1$ topology. The goal of this work is to get a positive answer of the second one in the $C^1$ topology. Generalizations of singular hyperbolicity in higher dimension have been proposed~\cite{MM,ZGW},
but it is not clear for us what should be the generalization of Conjecture~\ref{conj2}.
Since Newhouse phenomenon requires $C^2$-smoothness, an even stronger result could be
imagined, specific to the $C^1$-topology:
singular hyperbolicy may be $C^1$-dense in the space of three dimensional vector fields.

\section{Precise statements}

Let $M$ be a three-dimensional compact Riemannian manifold without boundary. A smooth vector field $X$ on $M$ generates a flow $\varphi_t$. A point $x$ is \emph{regular} if $X(x)\neq0$; otherwise, $x$ is \emph{singular}. The set of singularities plays a particular role and is denoted by ${\rm Sing}(X)$. The derivative ${\rm D}\varphi_t$ of $\varphi_t$ w.r.t. the space variable is called the \emph{tangent flow}.
\medskip

The dynamics of the flow is usually split in the following way.
For any $\varepsilon>0$, an \emph{$\varepsilon$-pseudo orbit from $x$ to $y$} is a sequence
$\{x_i\}_{i=0}^n$ such that $x_0=x$, $x_n=y$, $n\geq 1$ and such that
$d(\varphi_{t_i}(x_i),x_{i+1})<\varepsilon$ for any $i\in \{0,\dots,n-1\}$ and for some $\{t_i\}_{i=0}^{n-1}$ in $[1,2]$. Any two points $x,y\in M$ are said to be \emph{chain related} if for any $\varepsilon>0$, there are $\varepsilon$-pseudo orbits from $x$ to $y$ and from $y$ to $x$. If $x$ is chain related to itself, then $x$ is called a \emph{chain recurrent point}. The (compact invariant) set of chain recurrent points is called the \emph{chain recurrent set} of $X$.
To be chain related is an equivalence relation on the chain recurrent set;
the equivalence classes are called \emph{chain recurrence classes} of $X$.
They are compact, invariant (and pairewise disjoint).
\medskip

Let $\Lambda$ be an invariant compact set.
An invariant continuous splitting $T_\Lambda M=E\oplus F$
(where $E,F$ are non-trivial vector bundles) is \emph{dominated}
if there exist constants $C,\lambda>0$ such that
for any $t>0$, any $x\in\Lambda$ and for any unit vectors $u\in E(x)$ and $v\in F(x)$, we have
$\|{\rm D}\varphi_t(u)\|\le C{\rm e}^{-\lambda t}\|{\rm D}\varphi_t(v)\|$.

The set $\Lambda$ is \emph{hyperbolic}, if there exist a continuous invariant splitting $T_\Lambda M=E^s\oplus \left< X\right>\oplus E^u$ w.r.t. ${\rm D}\varphi_t$
and constants $C,\lambda>0$ such that for any $x\in\Lambda$ and $t>0$, one has $\|{D\varphi_t}|_{E^s(x)}\|\le C{\rm e}^{-\lambda t}$ and $\|{D\varphi_{-t}}|_{E^u(x)}\|\le C{\rm e}^{-\lambda t}$.

An attractor $\Lambda$ is \emph{singular hyperbolic} if there is a dominated splitting $T_\Lambda M=E^{ss}\oplus E^{cu}$ w.r.t. ${\rm D}\varphi_t$ and constants $C>0$ and $\lambda>0$
such that:
\begin{itemize}

\item Contraction: for any $t>0$ and any $x\in\Lambda$, $\|{\rm D}\varphi_t|_{E^{ss}(x)}\|\le C{\rm e}^{-\lambda t}$.

\item Area-expansion: for any $t>0$ and any $x\in\Lambda$, $|{\rm Det}\varphi_{-t}|_{E^{cu}(x)}|\le C{\rm e}^{-\lambda t}$.

\end{itemize}A transitive repeller is \emph{singular hyperbolic} if it is a singular hyperbolic
attractor for $-X$.
\medskip

We say that $X$ is \emph{singular hyperbolic} if the chain recurrent set of $X$ is the union of finitely pairewise disjoint invariant compact sets $\{\Lambda_i\}$ such that each $\Lambda_i$ is a hyperbolic set, a singular hyperbolic attractor, or a singular hyperbolic repeller.
$X$ is \emph{robustly} singular hyperbolic if any vector field $C^1$-close to $X$
is singular hyperbolic.
\medskip

We say that $X$ has a \emph{homoclinic tangency} if $X$ has a (non-singular) hyperbolic periodic orbit $\gamma$ such that the stable manifold of $\gamma$ and the unstable manifold of $\gamma$ have some non-transverse intersection.
\bigskip

We announce an answer to the above Palis conjecture~\ref{conj2} is:
\medskip

\begin{maintheorem}
In the $C^1$ topology, every three dimensional vector field can be accumulated by robustly singular hyperbolic vector fields, or by vector fields with homoclinic tangencies.
\end{maintheorem}
\medskip

For proving this theorem, we only need to prove the following (which was
already obtained by Arroyo and Rodriguez-Hertz \cite{ArR03} for non-singular chain-recurrence classes).
\medskip

\begin{maintheoremrestated}\label{Thm:semiglobal}
For $C^1$ generic vector field $X$ on a three-dimensional manifold which is far from homoclinic tangencies,
any chain-recurrence class is hyperbolic or is a singular hyperbolic attractor or repeller.
\end{maintheoremrestated}

\section{Dominated splittings for tangent \emph{vs} linear Poincar\'e flows}

The dynamics induces another linear flow above the set of regular points.
At any point $x\in M\setminus {\rm Sing}(X)$ we introduce the plane
$\cN_x=X(x)^\perp$ and define the normal bundle ${\cal N}=\coprod_{x\in M\setminus{\rm Sing}(X)}{\cal N}_x$. By orthogonal projection of the tangent flow ${\rm D}\varphi_t$,
one gets the \emph{linear Poincar\'e flow} $\psi_t$ on $\cN$.
\medskip

Let $C$ be any chain-recurrence class (which is not an isolated singularity)
of a $C^1$-generic vector fields far from homoclinic tangencies.
Using technics developed in the different works on robustly transitive sets~\cite{MPP04,LGW05,MoP03}
and Liao's estimation~\cite{Lia89},
Gan and Yang~\cite{GaY13} have proved the following properties (up to replace $X$ by $-X$):

\begin{enumerate}
\item 
Any non-isolated singularity $\sigma\in C$ is \emph{Lorenz-like}:
it has three real eigenvalues satisfying $\lambda_1<\lambda_2<0<-\lambda_2<\lambda_3$.

Moreover the (one-dimensional) strong stable manifold satisfies $W^{ss}(\sigma)\cap C=\{\sigma\}$.

\item The linear Poincar\'e flow on $C\setminus {\rm Sing(X)}$ has a dominated splitting.

\item  {\cite[Theorem C]{GaY13}.}  If  the tangent flow ${\rm D}\varphi_t$ on $\Lambda$ has a dominated splitting
and if $C$ contains a singularity, then $C$ is a singular hyperbolic attractor.

\end{enumerate}
\medskip

Thus our main theorem essentially reduces to compare the dominated splitting
of the linear Poincar\'e flow and of the tangent flow.
(Note that the following result goes beyond $C^1$-generic vector fields.)
\medskip

\begin{Theorem} \emph{Dominated splitting for the tangent flow. --}
\label{Thm:flowcontraction}
Consider any $C^3$ three-dimensional vector field $X$ and any compact invariant set $\Lambda$ with the following properties:

\begin{itemize}

\item Every periodic orbit in $\Lambda$ is a hyperbolic saddle.

\item Every singularity $\sigma\in \Lambda$ is Lorenz-like and $W^{ss}(\sigma)\cap \Lambda=\{\sigma\}$.

\item $\Lambda$ does not contain a minimal repeller whose dynamics is the suspension of
an irrational rotation of the circle.

\end{itemize}
Then the tangent flow ${\rm D}\varphi_t$ on $\Lambda$ has a dominated splitting
if and only if the linear Poincar\'e flow on $\Lambda\setminus {\rm Sing(X)}$ has a dominated splitting.

%
%
%
%

\end{Theorem}

%
%
%
%
%
%

\section{Uniform contraction for dominated fibered dynamics}

The existence of singularities introduces several difficulties:
the regular orbits may separate when passing the singularities,
which causes a lack of compactness and of uniformity.
These problems have been overcome in previous works in two ways:
\begin{itemize}
\item by blowing up the singular set and extending the flow (introduced by Li-Gan-Wen~\cite{LGW05}),
\item by rescaling the flow near the singular set (in Liao's work~\cite{Lia89} and in~\cite{GaY13}).
\end{itemize}
These ideas may be applied to several flows associated to the initial flow $\varphi_t$:
the flow itself, the sectional Poincar\'e flow, their tangent flows,...
\medskip

In this work we get a compactification of the rescaled sectional Poincar\'e flow
that we define now.
For any regular points $x$ and time $t_0\in \RR$,
the flow $\varphi_t$ induces a local holonomy map from a neighborhood of $x$ in
$\exp_x({\cal N}_x)$ to a neighborhood of $\varphi_t(x)$ in $\exp_{\varphi_t(x)}({\cal N}_{\varphi_t(x)})$. This induces a local map $P_{t_0}$ from ${\cal N}_{x}$ to ${\cal N}_{\varphi_t(x)}$
which preserves $0$ and gives a local flow $P_t$ on the bundle $\cN$
in a neighborhood of the $0$-section. It is called the \emph{sectional Poincar\'e flow}.
Its linearization at the $0$-section is the linear Poincar\'e flow $\psi_t$.

For these flows, we define the rescaled sectional Poincar\'e flow and the
rescaled linear Poincar\'e flow as:

$$P_t^*(v_x)=\|X(\varphi_t(x))\|^{-1}\cdot P_t\big(\|X(x)\|\cdot v_x\big) \text{ and }
\psi_t^*(v_x)=\frac{\|X(x)\|}{\|X(\varphi_t(x))\|}\cdot\psi_t(v_x).$$
Still $\psi_t^*$ can be seen as the linearization of $P_t^*$.
\medskip

\begin{Theorem} \emph{Compactification of the sectional Poincar\'e flow. --}
Assume that the flow $\varphi_t$ is $C^r$, $r\geq 2$,
and let $K$ be an invariant compact set whose singularities
${\rm Sing}(X)\cap K$ are hyperbolic.
Then, there exists a flow $\widehat \varphi_t$ on a compact metric space $\widehat K$
and a local flow $\widehat P^*_t$ on a linear bundle $\widehat \cN$ over
the dynamics of $\widehat \varphi_t$ on $\widehat K$,
which preserves the $0$-section, such that:
\begin{itemize}
\item There exists an injective continuous fiber-preserving map $p\colon  \cN\to \widehat \cN$
which conjugates the local flows $P^*_t$ and $\widehat P^*_t$.
\item $\widehat P^*_t$ is $C^{r-1}$-along fibers: it induces local $C^{r-1}$-diffeomorphisms
$(\widehat \cN_x,0) \to (\widehat \cN_{\widehat \varphi_t(x)},0)$ depending continuously
on $x$ for the $C^{r-1}$-topology.
\end{itemize}
\end{Theorem}
\bigskip

Theorem~\ref{Thm:flowcontraction}
can be translated for the flow $\widehat P^*_t$ with the following remarks:
\begin{itemize}
\item If the linear Poincar\'e flow $\psi_t$ has a dominated splitting,
then the $0$-section in $\widehat \cN$ has a dominated splitting
$T\widehat \cN=\widehat E\oplus \widehat F$ for the compactified rescaled flow $\widehat P^*_t$.
\item If moreover $\widehat E$ is uniformly contracted by the linearization of
$\widehat P^*_t$, then $K$ has a dominated splitting $T_KM=E\oplus F$ for the tangent
flow ${\rm D}\varphi_t$, with $\dim(E)=1$ and $\RR.X\subset F$.
\end{itemize}
\medskip

Theorem~\ref{Thm:flowcontraction} is thus a consequence of the theorem below for dominated local fibered flow $\widehat P^*_t$
on a bundle $\widehat \cN$ with $2$-dimensional fibers.
It requires an important \emph{compatibility assumption} which reflects the fact that $\widehat P^*_t$
has been obtained from the sectional Poincar\'e flow on a manifold:
\begin{enumerate}
\item There exists an open set $U\subset \widehat K$ and for close points $x,y\in U$
an identification map $\pi_{x,y}\colon \widehat \cN_x\to \widehat \cN_y$
which is ``compatible" with the local flow $\widehat P^*_t$.
\item Along pieces of orbit in $\widehat K\setminus U$, the bundle $\widehat E$ is uniformly contracted.
\end{enumerate}
\medskip

\begin{Theorem} \emph{Uniform contraction of one-dimensional bundle
for dominated local fibered flow. --}
Let $\widehat P^*_t$ be a local fibered flow which is $C^2$-along fibers,
on a bundle $\widehat \cN$ with $2$-dimensional fibers satisfying:
\begin{itemize}
\item the $0$-section is invariant and has a dominated splitting $T\widehat \cN=\widehat E\oplus \widehat F$;
\item the compatibility assumption holds;
\item $E$ is contracted along each periodic orbit of the base space $\widehat K$ of $\widehat \cN$;
\item the base space $\widehat K$ of $\widehat \cN$ does not contain any
invariant compact set which is a repeller and whose dynamics is conjugated to the suspension
of an irrational circle rotation.
\end{itemize}
Then $\widehat E$ is uniformly contracted by the linearization of the flow $\widehat P^*_t$.
\end{Theorem}
\bigskip

This theorem is a continuation of a sequence of works initiated by Ma\~n\'e~\cite{Man85}
and Pujals-Sambarino~\cite{PuS00}:
the later proves the hyperbolicity of dominated invariant compact set for surface diffeomorphisms.
Compared with the usual Pujals-Sambarino arguments, we meet the following difficulties:

\begin{itemize}

\item We work with continuous time dynamics. It is not enough to consider
the time-one map since the flow shears along the orbits.
Bounded time evolution produces small shear, but we need to consider long-term behaviors.
This difficulty already appears in Arroyo and Rodriguez-Hertz' result~\cite{ArR03}.
%

\item There is no extension of Pujals-Sambarino result
to general dominated fibered dynamics (as one can see considering
a product of a minimal base dynamics with the identity along fibers).
Our proof uses strongly the identification maps $\pi_{x,y}$.

\item Since the identifications $\pi_{x,y}$ are only defined on $U$,
we have to handle with the induced dynamics in $U$.
This makes us to consider ``induced hyperbolic returns'' (as in \cite{CrP11}).

\item The flow $\widehat P^*_t$ is locally defined and the global arguments of~\cite{PuS00}
have to be replaced by local ones.

\item We need to construct some Markovian boxes in a non-symmetric setting:
the bundles $\widehat E, \widehat F$ play different roles. We borrow some ideas from a recent work of Crovisier, Pujals and Sambarino \cite{CPS14}.

\end{itemize}


\begin{thebibliography}{1d}

\bibitem{ABS77}
V. Afra\u\i movi\v c, V. Bykov, and L. Silnikov, The origin and structure of the Lorenz attractor, {\it Dokl.
Akad. Nauk SSSR}, {\bf 234} (1977), 336--339.

\bibitem{ArR03}A. Arroyo and F. Rodriguez Hertz, Homoclinic bifurcations and uniform hyperbolicity for threedimensional
flows, {\it Ann. Inst. H. Poincar\'e-Anal. NonLin\'eaire}, {\bf20}(2003), 805--841.






\bibitem{CrP11} S. Crovisier and E. Pujals,
Essential hyperbolicity and homoclinic bifurcations: a dichotomy phenomenon/mechanism for
diffeomorphisms. {\it Preprint} ArXiv:1011.3836.

\bibitem{CPS14} S. Crovisier, M. Sambarino and E. Pujals, Hyperbolicity of the extremal bundles, {\it in preparation}.

%




\bibitem{GaY13} S. Gan and D. Yang, Morse-Smale systems and non-trivial horseshoes for three-dimensional singular flows. {\it Preprint} Arxiv:1302.0946v1.

\bibitem{Guc76}
J. Guckenheimer, A strange, strange attractor, {\it The Hopf bifurcation  theorems and its applications (Applied
Mathematical Series, {\bf 19})}, Springer-Verlag (1976), 368--381.


\bibitem{LGW05} M. Li, S. Gan and L. Wen, Robustly transitive singular sets via
approach of extended linear Poincar\'e flow, {\it Discrete Contin.
Dyn. Syst.}, {\bf 13 }(2005), 239--269.


\bibitem{Lia89} S. Liao, On $(\eta,d)$-contractible orbits of vector
fields, {\it Systems Science and Mathematical Sciences},
{\bf2} (1989), 193--227.


\bibitem{Man85} R. Ma{\~n}{\'e}, Hyperbolicity, sinks and measure in one-dimensional dynamics, {\it Comm. Math. Phys.}, {\bf 100} (1985), 495--524.

\bibitem{MM} R. Metzger and C. Morales,
Sectional-hyperbolic systems,
{\it Ergodic Theory Dynam. Systems}, {\bf 28} (2008), 1587--1597. 

\bibitem{MoP03} C. Morales and M. Pacifico,  A dichotomy for
three-dimensional vector fields, {\it Ergodic Theory and Dynamical
Systems}, {\bf23} (2003), 1575--1600.

\bibitem{MPP04}
C. Morales, M. Pacifico, and E. Pujals, Robust transitive singular sets for 3-flows are partially hyperbolic
attractors or repellers, {\it  Ann. of Math.}, {\bf160} (2004), 375--432.

\bibitem{New79} S. Newhouse,
The abundance of wild hyperbolic sets and nonsmooth stable sets for diffeomorphisms,
{\it Inst. Hautes \'Etudes Sci. Publ. Math.}, {\bf 50} (1979), 101--151. 


\bibitem{Pal91} J. Palis, Homoclinic bifurcations, sensitive-chaotic dynamics and strange attractors, {\it Dynamical systems and related topics (Nagoya, 1990)}, 
Adv. Ser. Dynam. Systems, \textbf{9}, {\it World Sci. Publ., River Edge, NJ}, (1991), 466-472, 

\bibitem{Pal00} J. Palis, A global view of dynamics and a conjecture of
the denseness of finitude of attractors, {\it Ast\'erisque}, {\bf
261} (2000), 335--347.

\bibitem{Pal05} J. Palis, A global perspective for non-conservative dynamics,
{\it Ann. Inst. H. Poincar\'e Anal. Non Lin\'eaire}, {\bf22} (2005),
485--507.






\bibitem{PuS00} E. Pujals and M. Sambarino, Homoclinic tangencies and
hyperbolicity for surface diffeomorphisms, {\it Annals of Math.},
{\bf 151} (2000), 961--1023.




%

\bibitem{ZGW} S. Zhu, S. Gan and L. Wen, Indices of singularities of robustly transitive sets,
{\it Discrete Contin. Dyn. Syst.}, {\bf 21} (2008),  945--957.

\end{thebibliography}
\end{document}